\documentclass[10pt]{amsart}
\usepackage{amsfonts}
\usepackage{amsmath}

\newcommand{\naturals}{\ensuremath{\mathbb{N}}}
\newcommand{\integers}{\ensuremath{\mathbb{Z}}}

\newcommand{\reals}{\ensuremath{\mathbb{R}}}
\newcommand{\complex}{\ensuremath{\mathbb{C}}}

 \newcommand{\isom}{\cong}
 
 \newcommand{\union}{\cup}

\newtheorem{theorem}{Theorem}[section]
\newtheorem{proposition}[theorem]{Proposition}
\newtheorem{lemma}[theorem]{Lemma}
\newtheorem{corollary}[theorem]{Corollary}

\newtheorem{definition}{Definition}[section]
\newtheorem{question}{Question}

\newtheorem{example}{Example}[section]

\title{Busemann Points of Infinite Graphs}

\author{Corran Webster \and Adam Winchester}

\address{Department of Mathematical Sciences, University of Nevada Las
Vegas, Las Vegas, NV 89154}

\email{cwebster@unlv.nevada.edu}

\date{28th July, 2003}



\begin{document}

\begin{abstract}
    We provide a geometric condition which determines whether or not
    every point on the metric boundary of a graph with the standard
    path metric is a Busemann point, that is it is the limit point of
    a geodesic ray.  We apply this and a related condition to
    investigate the structure of the metric boundary of Cayley graphs. 
    We show that every metric boundary point of the Cayley graph of a
    finitely generated Abelian group is a Busemann point, but groups
    such as the braid group and the discrete Heisenberg group have
    boundary points of the Cayley graph which are not Busemann points
    when equipped with their usual generators.
\end{abstract}

\maketitle

The metric compactification of a metric space was introduced by
Gromov~\cite{gromov:1}, but was little studied.  Recently, this
compactification has proven of use in the study of certain metrics on
state spaces of C*-algebras~\cite{rieffel:groupalgebras}.

Let $G$ be a countable discrete group, and $\complex_{c}(G)$ the
convolution algebra of functions with finite support.  If $\pi_{l}$ is
the usual *-representation of $\complex_{c}(G)$ on $\ell^{2}(G)$
coming from the left-regular representation, the reduced group
C*-algebra $C^{*}_{r}(G)$ is the closure of $\pi_{l}(\complex_{c}(G))$
in $B(\ell^{2}(G))$.  Given a length function $\ell$ on the group,
Connes~\cite{connes:compactmetric} considers as a ``Dirac'' operator
the unbounded operator $M_{\ell}$ on $\ell^{2}(G)$ given by
multiplication by $\ell$.  The commutators $[M_{\ell}, \pi_{l}(f)]$
are bounded for $f \in \complex_{c}(G)$, and thus define a seminorm
$L_{\ell}(f) = \|[M_{\ell}, \pi_{l}(f)]\|$ on this dense *-subalgebra
of $C^{*}_{r}(G)$.  This in turn defines a dual metric
\[
  \rho_{L_{\ell}}(\varphi, \psi) = \sup \{|\varphi(f) - \psi(f)|: f
  \in \complex_{c}(G), L_{\ell}(f) \le 1\}
\]
on the state space $S(C^{*}_{r}(G))$ of $C^{*}_{r}(G)$.

Rieffel~\cite{rieffel:groupalgebras} asks whether these seminorms are
in fact Lip-norms, as defined in his earlier
papers~\cite{rieffel:groupmetrics,rieffel:noncommmetrics}.  A seminorm
such as $L_{\ell}$ is a Lip-norm if the topology given by the dual
metric $\rho_{L_{\ell}}$ on the state space coincides with the weak-*
topology.  His principal examples have $\ell$ as a word-length on
$\integers^{d}$ given by a set of generators, or coming from embedding
$\integers^{d}$ in $\reals^{d}$ and obtaining a length function from a
norm, and Rieffel shows that in these cases the seminorms are in fact
Lip-norms.  His method relies on the fact the the group's action on
the boundary of the metric compactification is amenable (as studied by
Anantharaman-Delarouche~\cite{anantharaman:amenability}), as well as
finiteness conditions on the orbits of the action on the boundary.

There are a number of obstacles to this approach, however, not least
of which is understanding the boundary of the metric compactification
in concrete terms.  The easiest definition of the metric
compactification is as the primitive ideal space of a certain
subalgebra of the C*-algebra of continuous, bounded functions on $G$
via Gelfand's theorem, but Rieffel shows that one can find the points
on the boundary as limits of weakly-geodesic rays.  In some cases,
such as for finitely generated free groups, all the boundary points
occur as the limits of rays satisfying a stronger condition.  Rieffel
calls such points Busemann points, and they play a significant role in
his discussion of the action of $\integers^{d}$ on the metric
boundary.  Rieffel raises the question of when all points on the
boundary of the metric compactification are Busemann points.

In this paper, we look into this question in the setting of path
metrics on infinite graphs, with a particular interest in Cayley
graphs.  By considering conditions under which minimal paths between
triples of vertices eventually meet, we are able to give a geometric
condition which determines whether or not all metric boundary points
are Busemann points.  More precisely, every point on the metric
boundary is a Busemann point if and only if given any pair of
vertices, there are minimal paths from each to any distant third
vertex which eventually share a tail.

We then turn to look specifically at Cayley graphs of finitely
generated groups.  We provide an example which shows that even for
finitely presented groups, there are Cayley graphs which have boundary
points which are not Busemann points.  We give a second if and only if
condition for a Cayley graph to have boundary points which are not
Busemann points.  This condition is easier to check than the general
condition, and we apply it to a number of examples.  We show that
every metric boundary point is a Busemann point for finitely generated
groups $G \isom F_{k}/N$ where $N$ is also finitely generated, as well
as providing a geometric proof of a result of Develin~\cite{develin:cayleycompactification} that every metric boundary point 
of a finitely generated Abelian group is a Busemann point.  Finally, we use the conditions to
show that the Cayley graphs of groups such as the braid groups and the
discrete Heisenberg groups, when given standard sets of generators,
have boundary points which are not Busemann points.

This paper came from an undergraduate research project between the
authors.  The authors would like to thank Michelle Schultz for
organizing the undergraduate research seminar at UNLV, and Marc
Rieffel for his suggestion that this area would be a fruitful topic
for undergraduate research and for his helpful comments on preliminary
versions of this paper.

\section{The Metric Compactification and Busemann Points}

Following Rieffel~\cite{rieffel:groupalgebras}, let $(X, d)$ be a
complete, locally compact metric space.  Let $C_{b}(X)$ be the
commutative, unital C*-algebra of bounded, continuos functions on $X$
with the uniform norm $\|\cdot\|_{\infty}$, and $C_{\infty}(X)$ the
closed subalgebra of functions which vanish at infinity.  We define
functions $\varphi_{y,z}: X \to \reals$ by
\[
  \varphi_{y,z}(x) = d(x,y) - d(x,z).
\]
It is immediate from the triangle equality that
\[
  \| \varphi_{y,z} \|_{\infty} \le d(y,z),
\]
and so $\varphi_{y,z} \in C_{b}(X)$.

We let $\mathcal{G}(X,d)$ be the closed subalgebra generated by
$C_{\infty}(X)$, the constant functions, and the functions
$\{\varphi_{y,z}| y,z \in X\}$.  Then $\mathcal{G}(X,d)$ is a
commutative, unital C*-algebra, so Gelfand's theorem tells us that
$\mathcal{G}(X,d) \isom C_{b}(\overline{X}_{d})$, where $\overline{X}_{d}$ is
the maximum ideal space (or equivalently, the set of pure states) of
$\mathcal{G}(X,d)$.  $\overline{X}_{d}$ is a compact topological space,
containing $X$ as an open subset in a natural way, so we call
$\overline{X}_{d}$ the \emph{metric compactification} of $(X,d)$.  The set
$\overline{X}_{d} \setminus X$ can be naturally thought of as the boundary
at infinity of the compactification, so we will call the set
$\partial_{d}X = \overline{X}_{d} \setminus X$ the \emph{metric boundary}
of $X$.

Rieffel showed that if we fix some base point $z_{0}$, and define
$\varphi_{y} = \varphi_{z_{0},y}$, then $\mathcal{G}(X,d)$ is
generated by $C_{b}(X)$, the constant functions, and the functions
$\{\varphi_{y}| y \in X\}$, because $\varphi_{y,z} = \varphi_{z} -
\varphi_{y}$.  Note that this does not depend on the choice of
$z_{0}$.

A more concrete definition of this boundary is of interest.  To that
end, the following concepts are introduced:

\begin{definition}\label{defn:geodesicrays}
    Let $(X,d)$ be a metric space, and $T$ an unbounded subset of
    $\reals^{+}$ containing $0$, and let $\gamma: T \to X$.  We say
    that
    \begin{enumerate}
	\item $\gamma$ is a \emph{geodesic ray} if
	\[
	  d(\gamma(s), \gamma(t)) = |s - t|
	\]
	for all $s$, $t \in T$.
	
	\item $\gamma$ is an \emph{almost-geodesic ray} if for every
	$\varepsilon > 0$, there is an integer $N$ such that
	\[
	  |d(\gamma(t), \gamma(s)) + d(\gamma(s), \gamma(0)) - t| <
	  \varepsilon
	\]
	for all $t$, $s \in T$ with $t \ge s \ge N$.
	
	\item $\gamma$ is a \emph{weakly-geodesic ray} if for every $y
	\in X$ and every $\varepsilon > 0$, there is an integer $N$
	such that
	\begin{equation*}
	  |d(\gamma(t), \gamma(0)) - t| < \varepsilon
	\end{equation*}
	and
	\begin{equation*}
	  |d(\gamma(t), y) - d(\gamma(s), y) - (t - s)| < \varepsilon
	\end{equation*}
	for all $t$, $s \in T$ with $t$, $s \ge N$.
    \end{enumerate}
\end{definition}

It is immediate that every geodesic ray is an almost-geodesic ray. 
Rieffel showed that every almost-geodesic ray is a weakly-geodesic
ray.  The significance of weakly-geodesic rays is that their limits
are the points of the metric compactification in most cases.  Recall
that a metric is \emph{proper} if every closed ball of finite radius
is compact.

\begin{theorem}[Rieffel]\label{thm:metricboundary}
    Let $(X,d)$ be a complete, locally compact metric space, and let
    $\gamma: T \to X$ be a weakly geodesic ray in $X$.  Then
    \[
      \lim_{t \to \infty} f(\gamma(t))
    \]
    exists for every $f \in \mathcal{G}(X,d)$, and defines an element
    of $\partial_{d}X$.  Conversely, if $d$ is proper and if $(X,d)$
    has a countable base, then every point of $\partial_{d}X$ is
    determined as above by a weakly-geodesic ray.
\end{theorem}

We will not reproduce the entire proof of the theorem, but the
construction involved in the last part will be of use to us later, so
we will reproduce that here.  The proof requires one additional result
from Rieffel's paper:

\begin{proposition}[Rieffel]\label{prop:countablebase}
    Let $(X,d)$ be a complete locally compact metric space.  If the
    topology of $X$ has a countable base, then so do the topologies of
    $\overline{X}_{d}$ and $\partial_{d}X$.
\end{proposition}

With this in hand we can prove the last part of the theorem:

\begin{proof}[Proof (Theorem~\ref{thm:metricboundary})]
    Let $\omega \in \partial_{d}X$. 
    Proposition~\ref{prop:countablebase} tells us that $\overline{X}_{d}$
    has a countable base, so we can find a sequence $w_{n} \in X$
    which converges to $\omega$ in $\overline{X}_{d}$.  Since $\omega
    \notin X$, and $d$ is proper, $w_{n}$ is unbounded.  So we can
    find a subsequence $w_{n_{k}}$ (with $w_{n_{0}} = w_{0}$) so that
    if $k > l$, then $d(w_{n_{k}}, w_{0}) > d(w_{n_{l}}, w_{0})$.  Let
    $T = \{ d(w_{n_{k}}, w_{0}) : k = 0, 1, \ldots \}$, and define
    $\gamma: T \to X$ by letting $\gamma(t) = w_{n_{k}}$ where $t =
    d(w_{n_{k}}, w_{0})$.  Clearly
    \[
      \lim_{t \to \infty} \gamma(t) = \omega,
    \]
    so we need only show that $\gamma$ is weakly geodesic.
    
    By construction, $d(\gamma(t),\gamma(0)) = t$, so $\gamma$
    satisfies the first condition of
    Definition~\ref{defn:geodesicrays}.3 for all $\varepsilon > 0$. 
    Use $\gamma(0)$ as the base point for functions $\varphi_{y}$ for
    $y \in X$.  Given any one of these functions we know that
    \[
      \lim_{t \to \infty} \varphi_{y}(\gamma(t)) =
      \varphi_{y}(\omega),
    \]
    and so given any $\varepsilon > 0$, there is some $N$ such that
    for all $s$, $t \in T$ with $s$, $t \ge N$ then
    $|\varphi_{y}(\gamma(s)) - \varphi_{y}(\gamma(t))| < \varepsilon$. 
    But then
    \begin{align*}
      |d(\gamma(t), y) - d(\gamma(s), y) - (t - s)| & = |d(\gamma(t),
      y) - d(\gamma(s), y) \\
      &\qquad - d(\gamma(t),\gamma(0)) + d(\gamma(s),\gamma(0)))|\\
      &= |\varphi_{y}(\gamma(t)) - \varphi_{y}(\gamma(s))| <
      \varepsilon,
    \end{align*}
    so the second condition for weakly-geodesic rays is satisfied.
\end{proof}

Note that this theorem and the above construction mean that given any
weakly-geodesic ray we can find a weakly-geodesic ray which has the
same limit in the metric boundary, but for which $d(\gamma(t),
\gamma(0)) = t$, and the ray comes from a sequence of points.

Rieffel defines any point $\partial_{d}X$ which is the limit of an
almost-geodesic ray to be a \emph{Busemann point}, and poses the
following question:

\begin{question}\label{quest:rieffel}
  Given a metric space $(X,d)$, is every point of $\partial_{d}X$ a
  Busemann point?
\end{question}

Rieffel's interest in the metric compactification came from looking at
metrics on infinite discrete groups.  The most natural metrics on
discrete groups are those which come from the standard graph metric on
a Cayley graph of the group.  So a natural class of metric spaces to
investigate are graph metrics.

\section{Graph Metrics}

If $\Gamma = (V, E)$ is a connected graph with vertices $V$ and edges
$E$, then the standard metric $d$ on $V$ is defined by letting
$d(x,y)$ be the minimum length of a path from $x$ to $y$.  Given
vertices $x$ and $y$, we will use the notation $[x,y]$ for an
unspecified, but fixed, minimal path from $x$ to $y$.

It is immediate that this metric gives $V$ the discrete topology, so
every function on $V$ is continuous, $(V,d)$ is automatically
complete, and is locally compact.  Furthermore $(V, d)$ has a
countable base if and only if $V$ is countable, and $d$ is proper if
and only if the closed ball of radius 1 about every vertex is a finite
set, or equivalently, if and only if every vertex has finite valence.

These restrictions on the graph are not very onerous at all.  In
particular most interesting Cayley graphs satisfy these restrictions:

\begin{example}
    Let $G$ be a finitely generated group with generating set $S =
    S^{-1}$.  Then the Cayley graph of $\Gamma_{G} = (G, E)$
    corresponding to this generating set is a connected graph, with
    $G$ countable, and every vertex has valence at most $|S|$.  Hence
    $(G, d)$ is a complete, locally compact, proper metric space with
    a countable base.
\end{example}

Because the metric takes on only integer values, we can make certain
simplifying assumptions.  We first note that the functions
$\varphi_{v}$ can only take integer values, and hence if $\gamma$ is a
weakly-geodesic ray which converges to some point $\omega \in
\partial_{d}V$,
\[
   \varphi_{v}(\omega) = \lim_{t \to \infty} \varphi_{v}(\gamma(t))
   \in \integers.
\]
Moreover, if $\gamma'$ is another weakly-geodesic ray it converges to
$\omega$ if and only if for every $v \in V$, $\varphi_{v}(\gamma'(t))
= \varphi_{v}(\gamma(t))$ eventually.

Also it is easy to see that if $\gamma: T \to V$ is a geodesic ray,
then $T \subseteq \integers^{+}$, and we can in fact extend the domain
to all of $\integers^{+}$: if $T = \{t_{0} = 0, t_{1}, t_{2}, \ldots
\}$, we simply find minimal paths $[0,t_{1}]$, $[t_{1},t_{2}]$,
$[t_{2},t_{3}]$,\ldots and let $\gamma(n)$ be the $(n-t_{k-1})$th
point on the $k$th path for $t_{k-1} < n < t_{k}$.

Furthermore, it turns out that for these metrics, every Busemann point
is in fact the limit of a geodesic ray.

\begin{lemma}
  Let $\Gamma = (V,E)$ be a connected graph, and $d$ the usual graph
  metric.  Then if $\omega \in \partial_{d}V$ is a Busemann point,
  there is a geodesic ray $\gamma: T \to V$ which converges to
  $\omega$.
\end{lemma}

\begin{proof}
  Let $\gamma': T' \to V$ be an almost-geodesic ray which converges to
  $\omega$.  We can find an $N$ such that for all $s$, $t \in T'$ with
  $t \ge s \ge N$, we have
  \[
    |d(\gamma'(t),\gamma'(s)) + d(\gamma'(s),\gamma'(0)) - t| < 1/3.
  \]
  In particular, $|d(\gamma'(t),\gamma'(0)) - t| < 1/3$.  We let
  $t_{0} = t_{0}' = 0$, and find $t_{n} \in \naturals$ and $t_{n}' \in
  T'$ so that $t_{n} > t_{n-1}$, $t_{n}' \ge N$ and
  $d(\gamma'(t_{n}'),\gamma'(0)) = t_{n}$.  Let $T = \{t_{0}, t_{1},
  t_{2}, \ldots \}$, and define $\gamma: T \to V$ by $\gamma(t_{n}) =
  \gamma(t_{n}')$.  This implies that $|t_{n} - t_{n}'| =
  |d(\gamma'(t_{n}),\gamma'(0)) - t_{n}'| < 1/3$, and so $t_{n} \ge
  t_{n}' - 1/3$.
  
  Now for any $t_{n}$, $t_{m} \in T$, with $n \ge m$, we have that
  \begin{align*}
    |d(\gamma(t_{n}),\gamma(t_{m})) - (t_{n} - t_{m})| &\le
    |d(\gamma(t_{n}),\gamma(t_{m})) + t_{m} - t_{n}'| + 1/3\\
    &\le |d(\gamma'(t_{n}'),\gamma'(t_{m}')) +
    d(\gamma'(t_{m}'),\gamma'(0)) - t_{n}'| + 1/3\\
    &\le 2/3 < 1.
  \end{align*}
  But both $d(\gamma(t_{n}),\gamma(t_{m}))$ and $(t_{n} - t_{m})$ are
  integers, so
  \[
    d(\gamma(t_{n}),\gamma(t_{m})) = |t_{n} - t_{m}|
  \]
  and we conclude that $\gamma$ is a geodesic ray.
\end{proof}

In particular, this lemma tells us that if we wish to show that a
point on the boundary is not a Busemann point, it is sufficient to
show that it is not the limit of any geodesic ray.  The following
simple example uses this fact to show that we cannot hope to answer
Question~\ref{quest:rieffel} in the affirmative for general graphs.

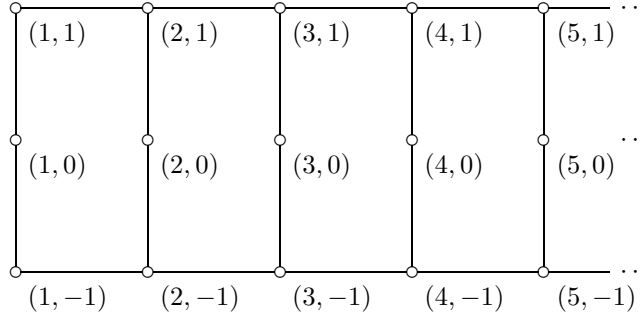
\begin{figure}
    \begin{picture}(260,110)(-10,-60) \multiput(-3,-50)(50,0){4}{
    \begin{picture}(50,100)(0,-50) \put(0,0){\circle{4}}
    \put(0,50){\circle{4}} \put(0,-50){\circle{4}}
    \put(0,2){\line(0,1){46}} \put(0,-2){\line(0,-1){46}}
    \put(2,50){\line(1,0){46}} \put(2,-50){\line(1,0){46}}
	\end{picture}}

	\put(200,0){\circle{4}} \put(200,50){\circle{4}}
	\put(200,-50){\circle{4}} \put(200,2){\line(0,1){46}}
	\put(200,-2){\line(0,-1){46}} \put(202,50){\line(1,0){23}}
	\put(202,-50){\line(1,0){23}}

	\put(5,45){\makebox(0,0)[tl]{$(1,1)$}}
	\put(55,45){\makebox(0,0)[tl]{$(2,1)$}}
	\put(105,45){\makebox(0,0)[tl]{$(3,1)$}}
	\put(155,45){\makebox(0,0)[tl]{$(4,1)$}}
	\put(205,45){\makebox(0,0)[tl]{$(5,1)$}}
	
	\put(5,-5){\makebox(0,0)[tl]{$(1,0)$}}
	\put(55,-5){\makebox(0,0)[tl]{$(2,0)$}}
	\put(105,-5){\makebox(0,0)[tl]{$(3,0)$}}
	\put(155,-5){\makebox(0,0)[tl]{$(4,0)$}}
	\put(205,-5){\makebox(0,0)[tl]{$(5,0)$}}
	
	\put(5,-55){\makebox(0,0)[tl]{$(1,-1)$}}
	\put(55,-55){\makebox(0,0)[tl]{$(2,-1)$}}
	\put(105,-55){\makebox(0,0)[tl]{$(3,-1)$}}
	\put(155,-55){\makebox(0,0)[tl]{$(4,-1)$}}
	\put(205,-55){\makebox(0,0)[tl]{$(5,-1)$}}
	
	\put(229,50){\makebox(0,0)[l]{$\cdots$}}
	\put(229,0){\makebox(0,0)[l]{$\cdots$}}
	\put(229,-50){\makebox(0,0)[l]{$\cdots$}}

    \end{picture}
    \caption{\label{fig:nonbusemannpoints}The graph of
    Example~\ref{eg:nonbusemannpoints}}
\end{figure}

\begin{example}\label{eg:nonbusemannpoints}
    Let $ \Gamma = (V,E)$ be the graph where $V = \naturals \times
    \{-1,0,1\}$, and there are edges joining
    \begin{itemize}
	\item $(k,1)$ and $(k+1,1)$ for $k \in \naturals$,
	
	\item $(k,-1)$ and $(k+1,-1)$ for $k \in \naturals$,
	
	\item $(k,1)$ and $(k,0)$ for $k \in \naturals$,
	
	\item $(k,-1)$ and $(k,0)$ for $k \in \naturals$.
    \end{itemize}
   
    This graph is illustrated in Figure~\ref{fig:nonbusemannpoints}.
    
    Simple calculations show that with $(1,0)$ as the base point and
    $l > k$, we have:
    \begin{align*}
      \varphi_{(k,1)}(x) &= \begin{cases} k & \text{for } x = (l,1) \\
	k & \text{for }x = (l,0) \\
	k-2 & \text{for }x = (l,-1) \end{cases}, \\
      \varphi_{(k,0)}(x) &= k-1, \\
      \varphi_{(k,-1)}(x) &= \begin{cases} k-2 & \text{for }x = (l,1)
      \\
	k & \text{for }x = (l,0) \\
	k & \text{for }x = (l,-1) \end{cases}.
    \end{align*}
    Hence the functions in $\mathcal{G}(V, d)$ which do not vanish at
    infinity are, for each $j \in \{1,0,-1\}$, eventually constant for
    $(l,j)$ as $l$ goes to infinity, although the constants are
    different for different $j$.
    
    Hence the metric boundary consists of just 3 points given by the
    following three weakly-geodesic rays:
    \begin{enumerate}
	\item $\gamma_{1}: \naturals \to \Gamma$, where $\gamma_{1}(n)
	= (n,1)$,
    
	\item $\gamma_{0}: \naturals \to \Gamma$, where $\gamma_{0}(0)
	= (1,1)$, and $\gamma_{0}(n) = (n,0)$,

	\item $\gamma_{-1}: \naturals \to \Gamma$, where
	$\gamma_{-1}(n) = (n,-1)$.
    \end{enumerate}
    Notice that $\gamma_{1}$ and $\gamma_{-1}$ are in fact geodesic
    rays; but $\gamma_{0}$ is not.
    
    In fact there is no geodesic ray which gives the same limit at
    infinity as $\gamma_{0}$, since points on any geodesic ray must
    eventually be either of the form $(k,1)$ or $(k,-1)$.  Hence
    $\lim_{n \to \infty} \gamma_{0}(n)$ is not a Busemann point.
\end{example}

It is perhaps of interest that all triangles this example are $2$-slim
(as in Alonso~\cite{alonso:wordhyperbolicgroups}, for example), so this is a
hyperbolic metric space.

Note that in the above example, there is no commonality in the minimal
paths from $(1,1)$ to $\gamma_{0}(n)$ and the minimal paths from
$(1,-1)$ to $\gamma_{0}(n)$.  Our key result is that this sort of
situation precisely characterizes when there are metric boundary
points which are not Busemann points.

Before proceeding, we make the following definition:

\begin{definition}
    Let $\Gamma = (V, E)$ be a connected graph and $d$ the graph
    metric.  The \emph{perimeter} of a triple of vertices $\{a,b,c\}$
    is $d(a,b) + d(b,c) + d(c,a)$.
\end{definition}

In the example, we have a pair of vertices $a$ and $b$ for which we
can find a sequence of vertices $c_{n}$ such that the perimeters of
the triples $\{a,b,c_{n}\}$ get arbitrarily large, and there are no
minimal paths from $a$ and $b$ to $c_{n}$ which share a tail.  This
following theorem tells us that if there is such a pair of points,
then this guarantees the existence of boundary points which are not
Busemann points.

\begin{theorem}\label{thm:nonbusemann}
    Let $\Gamma = (V,E)$ be a connected graph where $V$ is countable
    and every vertex has finite valence, such that there is a pair of
    vertices $a$, $b$ such that for every $n \in \naturals$, there is
    a vertex $c_{n} \in V$ such that the triple $\{a,b,c\}$ has
    perimeter $m \ge n$, and no minimal paths $[a,c]$ and $[b,c]$
    share a tail.  Then there is a point in $\partial_{d}V$ which is
    not a Busemann point.
\end{theorem}
\begin{proof}
    Without loss of generality, we may assume that $a$ is the base
    point for the functions $\varphi_{v}$.
    
    The vertices $c_{n}$ are a sequence in $\overline{X}_{d}$, which is
    compact, so there is a convergent subsequence $c_{n_{k}}$.  Since
    the perimeters are getting bigger, the points $c_{n_{k}}$ must
    head to infinity, and hence the limit point of $c_{n_{k}}$ is an
    element $\omega$ of $\partial_{d}V$.  By the same argument as
    Theorem~\ref{thm:metricboundary}, we can find a weakly-geodesic
    ray $\gamma: T \to V$, corresponding to a subsequence of
    $c_{n_{m}}$, which converges to that point on the boundary. 
    Without loss of generality, we may assume that $\gamma(0) = a$.
    
    Now assume that $\omega$ is a Busemann point, so we can find a
    geodesic ray $\gamma': \naturals \to V$ with $\gamma'(0) = a$
    which converges to $\omega$.  For all $v \in V$, we must therefore
    have
    \[
      \lim_{t \to \infty} \varphi_{v}(\gamma(t)) = \lim_{t \to \infty}
      \varphi_{v}(\gamma'(t)) = \varphi_{v}(\omega).
    \]
    But since $\varphi_{v}(\omega)$ takes on integer values only, for
    each $v$ we can find some $N_{v}$, such that
    \[
      \varphi_{v}(\gamma(t)) = \varphi_{v}(\gamma'(t)) =
      \varphi_{v}(\omega)
    \]
    for all $t \in T$ with $t \ge N_{v}$.
    
    In particular, for $v = b$, let $y = \varphi_{b}(\omega)$.  Then
    \[
      d(\gamma'(t), b) = d(\gamma'(t), a) - \varphi_{b}(\gamma'(t)) = t - y
    \]
    for all $t \ge N_{b}$.  But also
    \[
      d(\gamma(t), b) = d(\gamma(t), a) - \varphi_{b}(\gamma(t)) = t - y
    \]
    for all $t \ge N_{b}$.
    
    Now, fix a particular $s \ge N_{b}$ and consider $v = \gamma'(s)$. 
    Then we have
    \begin{align*}
      \varphi_{\gamma'(s)}(\gamma(t)) &=
      \varphi_{\gamma'(s)}(\gamma'(t)) \\
      &= d(\gamma'(t),a) - d(\gamma'(t), \gamma'(s)) = t - (t-s) = s
    \end{align*}
    for $t \in T$ with $t \ge s$ and $t \ge N_{\gamma'(s)}$.  But this
    means that for such $t$,
    \[
      d(\gamma(t), \gamma'(s)) = d(\gamma(t),a) -
      \varphi_{\gamma'(s)}(\gamma(t)) = t - s,
    \]
    and so
    \[
      d(\gamma(t), \gamma'(s)) + d(\gamma'(s), a) = (t - s) + s =
      d(\gamma(t), a).
    \]
    In other words, there is a minimal path from $a$ to $\gamma(t)$
    which goes through $\gamma'(s)$.
    
    But also
    \[
      d(\gamma(t), \gamma'(s)) + d(\gamma'(s), b) = (t - s) + (s - y)
      = t - y = d(\gamma(t), b).
    \]
    So there is a minimal path from $b$ to $\gamma(t)$ which goes
    through $\gamma'(s)$.  So we have constructed two minimal paths
    which contradict our assumption for some $c_{n} = \gamma(t)$.
\end{proof}

We now show that this condition is sharp: if the perimeter of bad
triples like these is bounded for every pair of points, we can
guarantee that all boundary points are Busemann points.

\begin{theorem}\label{thm:busemann}
    Let $\Gamma = (V,E)$ be a connected graph with $V$ countable and
    each vertex has finite valence, and where for each pair of
    vertices $a$, $b$ there is some number $M_{a,b}$ such that if $c$
    is any vertex for which no minimal path from $a$ to $c$ shares a
    tail with a minimal path from $b$ to $c$, then the perimeter of
    $\{a,b,c\}$ is less than $M_{a,b}$.  Then every point on the
    metric boundary is a Busemann point.
\end{theorem}

\begin{proof}
    Let $\omega \in \partial_{d}V$, and let $\gamma: T \to V$ be a
    weakly geodesic ray which converges to $\omega$.  Without any loss
    of generality (using the construction from
    Theorem~\ref{thm:metricboundary}) we can assume that $T \subseteq
    \naturals$ and that $d(\gamma(t),\gamma(0)) = t$.  We seek a
    geodesic ray which converges to $\omega$.
    
    Since $V$ is countable, let $\{v_{k}: k \in \naturals\}$ be an
    enumeration of $V$.  Without loss of generality, we may assume
    that $v_{0} = \gamma(0)$.
    
    For each $n$ we will inductively find numbers $m_{n}$, vertices
    $w_{m_{n}}$ and subsequences $T_{n}$ of $T$ with the following
    properties:
    \begin{enumerate}
	\item $d(w_{m_{n}}, \gamma(0)) = m_{n}$.
	
	\item for all $t \in T_{n}$, and all $v_{k}$ for $k \le n$,
	there exists a minimal path from $v_{k}$ to $\gamma(t)$ which
	passes through $w_{m_{n}}$.
	
	\item if $n \ge 1$, then for all $t \in T_{n}$, there is a
	minimal path from $w_{m_{n-1}}$ to $\gamma(t)$ which passes
	through $w_{m_{n}}$.
    \end{enumerate}
    
    Let $m_{0} = 0$, $w_{0} = \gamma(0)$ and $T_{0} = T$.  All
    conditions are trivially satisfied in this case.
    
    Given $m_{n}$, $w_{m_{n}}$ and $T_{n}$, let $l =
    d(\gamma(0),v_{n+1})$, let
    \[
      M' > (M_{w_{m_{n}},v_{n+1}} - d(v_{n+1}, w_{m_{n}}))/2
    \]
    be a whole number, and let $m_{n+1} = m_{n} + M' + l$.
    
    For each $t \in T_{n}$, with $t > m_{n+1}$ we have that
    $d(v_{n+1}, w_{m_{n}}) \le m_{n} + l$ by the triangle inequality. 
    Also $d(v_{n+1}, \gamma(t)) \ge m_{n+1} - l$ and $d(w_{m_{n}},
    \gamma(t)) \ge m_{n+1} - m_{n}$, so
    \[
      d(v_{n+1}, \gamma(t)) + d(w_{m_{n}}, \gamma(t)) \ge 2m_{n+1} -
      m_{n} - l = 2M' + m_{n} + l > M_{w_{m_{n}},v_{n+1}}.
    \]
    So by the hypotheses there is some vertex $z_{t}$ where minimal
    paths from $v_{n+1}$ and $w_{m_{n}}$ to $\gamma(t)$ join, and we
    can assume that $d(z_{t},\gamma(0)) \le m_{n+1}$, since otherwise
    $\{w_{m_{n}},v_{n+1},z_{t}\}$ is a triple with perimeter greater
    than $M_{w_{m_{n}},v_{n+1}}$, and we can find a closer point where
    the minimal paths join.
    
    By following the minimal paths out to distance $m_{n+1}$ from
    $\gamma(0)$, we can have $d(z_{t},\gamma(0)) = m_{n+1}$ without
    loss of generality.  In fact, $z_{t}$ must be on a minimal path
    from every element of $\{v_{0}, v_{1}, \ldots, v_{n+1}\} \cup
    \{w_{m_{0}}, w_{m_{1}}, \ldots, w_{m_{n}}\}$ to $\gamma(t)$, since
    we can find a minimal path from points $v_{k}$ or $w_{m_{k}}$ to
    $\gamma(t)$ which passes through $w_{m_{n}}$, so we simply replace
    the tail of that path with the path from $w_{m_{n}}$ to
    $\gamma(t)$ which passes through $z_{t}$.
    
    Now these points $z_{t}$ may be different for different $t \in
    T_{n}$, but since each of these points lies in
    \[
      S_{m_{n+1}}(\gamma(0)) = \{ v \in V: d(\gamma(0),v) = m_{n+1}\},
    \]
    and this is a finite set since $d$ is proper, there must be at
    least one point $w_{m_{n+1}} \in S_{m_{n+1}}(\gamma(0))$ such that
    $w_{m_{n+1}} = z_{t}$ for infinitely many $t$.  Let $T_{n+1} = \{t
    : w_{m_{n+1}} = z_{t}\}$.
    
    By construction, $m_{n+1}$, $w_{m_{n+1}}$ and $T_{n+1}$ satisfy
    all three conditions.
    
    We claim that $\gamma': T' \to V$, where $T' = \{m_{n}: n \in
    \naturals\}$ and $\gamma'(m_{n}) = w_{m_{n}}$ is a geodesic ray. 
    This follows immediately from that fact that, by construction,
    given any $s$, $t \in T'$, with $t > s$, there is a minimal path
    from $\gamma(0)$ to $w_{t}$ which passes through $w_{s}$, and so
    \[
      d(w_{t},w_{s}) = d(w_{t},\gamma(0)) - d(w_{s},\gamma(0)) = t -
      s.
    \]
    
    Finally, we claim that this geodesic ray converges to $\omega$. 
    Given any $v_{n} \in V$, we know that for $t > m_{n}$, there is a
    minimal path from $v$ to $w_{t}$ which passes through $w_{m_{n}}$,
    and a minimal path from $\gamma(0)$ to $w_{t}$ which passes
    through $w_{m_{n}}$, and so
    \begin{align*}
      \varphi_{v_{n}}(w_{t}) &= d(w_{t},\gamma(0)) - d(w_{t},v_{n}) \\
      &= d(w_{t},w_{m_{n}}) + d(w_{m_{n}},\gamma(0)) -
      (d(w_{t},w_{m_{n}}) + d(w_{m_{n}},v_{n})) \\
      & = d(w_{m_{n}},\gamma(0)) - d(w_{m_{n}},v_{n}).
    \end{align*}
    Hence
    \[
      \lim_{t \to \infty} \varphi_{v_{n}}(w_{t}) =
      d(w_{m_{n}},\gamma(0)) - d(w_{m_{n}},v_{n}).
    \]
    
    On the other hand, given $t \in T_{n}$, there is a minimal path
    from $\gamma(0)$ to $\gamma(t)$ which passes through $w_{m_{n}}$
    and a minimal path from $v_{n}$ to $\gamma(t)$ which passes
    through $w_{m_{n}}$, and so
    \begin{align*}
      \varphi_{v_{n}}(\gamma(t)) &= d(\gamma(t),\gamma(0)) -
      d(\gamma(t),v_{n}) \\
      &= d(\gamma(t),w_{m_{n}}) + d(w_{m_{n}},\gamma(0)) -
      (d(\gamma(t),w_{m_{n}}) + d(w_{m_{n}},v_{n})) \\
      & = d(w_{m_{n}},\gamma(0)) - d(w_{m_{n}},v_{n}).
    \end{align*}
    Since $T_{n}$ is a subsequence of $T$,
    \[
      \varphi_{v_{n}}(\omega) = \lim_{t \in T_{n}}
      \varphi_{v_{n}}(\gamma(t)) = d(w_{m_{n}},\gamma(0)) -
      d(w_{m_{n}},v_{n}),
    \]
    and so $\varphi_{v_{n}}$ does not separate the limits of the two
    sequences.
    
    Since this happens for every $v_{n}$, we have that
    \[
      \omega = \lim_{n \to \infty} w_{m_{n}}.
    \]
    Hence $\omega$ is a Busemann point.
\end{proof}

It is perhaps worth noting that the proof of these result extends to
the slightly more general case of weighted graph metrics with edge
weights in $\naturals$, or in fact, $\lambda \naturals$ for any
$\lambda > 0$.

\begin{example}\label{eg:tree}
    Let $\Gamma = (V,E)$ be a tree.  Since the only way that the
    unique minimal paths $[a,c]$, and $[b,c]$ can fail to share a tail
    is if $c$ is on the unique minimal path $[a,b]$, we conclude that
    if $d(a,b) = n$, then the minimal paths must share a tail if the
    perimeter of the triple $\{a,b,c\}$ is greater than $2n$.  Hence
    every metric boundary point of a tree is a Busemann point.
\end{example}

\begin{example}\label{eg:lattice}
    Consider the integer lattice $V = \integers^{d} \subset
    \reals^{d}$ as the set of vertices, with edges joining vertices
    which differ by $\pm e_{k}$, where $e_{k}$ is a standard basis
    vector in $\reals^{d}$, so $\Gamma = (V, E)$ is the Cayley graph
    of $\integers^{d}$ with the standard generators.  Here there are
    many possible minimal paths between points, in general.  Given a
    fixed pair of points $a$ and $b$, and any point $c$, the only way
    that there can be no minimal path $[a,c]$ which shares a tail with
    a minimal path $[b,c]$ is if $c$ lies on some minimal path
    $[a,b]$.  So once again these minimal paths must share a tail if
    the perimeter of $\{a, b, c\}$ is greater than $2n$.  Hence every
    metric boundary point such a lattice is a Busemann point.
\end{example}

An obvious question arises concerning the way in which boundaries of
related graphs may be related.  For example, two metric spaces
$(X_{1},d_{1})$ and $(X_{2},d_{2})$ are Lipschitz equivalent, if there
is a bijection $T: X_{1} \to X_{2}$ and constants $\lambda_{1}$,
$\lambda_{2} > 0$ such that
\[
  \frac{d_{2}(T(a), T(b))}{d_{1}(a,b)} \le \lambda_{1} \qquad \text{and}
  \qquad \frac{d_{1}(a,b)}{d_{2}(T(a), T(b))} \le \lambda_{2}
\]
for all $a$ and $b \in X_{1}$, with $a \ne b$.  It is known from
Rieffel's work that Lipschitz equivalent metric spaces may have
different metric boundaries (even for graph metrics), and as the
following example shows, even if the metric boundaries are naturally
homeomorphic, which points are Busemann points may vary.

\begin{example}\label{eg:lipschitzequivalent}
    Let $\Gamma_{1}$ be the graph of
    Example~\ref{eg:nonbusemannpoints}, and let $\Gamma_{2}$ be the
    same graph, but with additional edges from $(k,0)$ to $(k+1,0)$. 
    The metric boundary of this second graph is homeomorphic to the
    metric boundary of the first graph and the boundary points are the
    limits of corresponding weakly-geodesic rays, but every metric
    boundary point of $\Gamma_{2}$ is a Busemann point.  However, the
    identity map on the vertices gives a Lipschitz equivalence between
    the two metric spaces, with
    \[
      \frac{d_{2}(a,b)}{d_{1}(a,b)} \le 1 \qquad \text{and} \qquad
      \frac{d_{1}(a,b)}{d_{2}(a,b)} \le 3.
    \]
\end{example}

\section{Cayley Graphs}

Our primary motivation is in the Cayley graphs of groups.  Recall that
$G = \langle S | R \rangle$ is a \emph{presentation} of a group if $S$
is a set of symbols, $R$ is a set of reduced words in $S$, and if $G$
is isomorphic to the group of equivalence classes of words in $S
\union S^{-1}$ where two words are equivalent if you can change one
word to another by adding and removing subwords of the form $ss^{-1}$
for $s \in S \union S^{-1}$, and subwords of the form $r$ or $r^{-1}$
for $r \in R$.  If this is the case, then $G \equiv F_{|S \union
S^{-1}|}/N$, where $F_{k}$ is the free group on $k$ generators, and
$N$ is the normal subgroup of the free group generated by $R$.  A
presentation is \emph{finite} if both $S$ and $R$ are finite sets.

Even in the case of Cayley graphs of finitely presented groups, we can
have metric boundary points which are not Busemann points.

\begin{example}\label{eg:nonbusemangroup}
    Consider the group $G = \langle a, b, c, d | aba^{-1}dcd^{-1}
    \rangle$.  Then the only minimal path from $a$ to $ab^{n}a^{-1}$
    is $b^{n}a^{-1}$, and the only minimal path from $d$ to
    $ab^{n}a^{-1} = dc^{n}d^{-1}$ is $c^{n}d^{-1}$.
    
    So $\{a, ab^{n}a^{-1}, d\}$ is a triple with perimeter $2n+2$, and
    so we have a metric boundary point which is not a Busemann point
    by Theorem~\ref{thm:nonbusemann}.
\end{example}

On the other hand, the Cayley graph of a finitely generated free group
with the standard generating set is a tree, so by
Example~\ref{eg:tree} we have that every metric boundary point is a
Busemann point.  Similarly the Cayley graph of the free Abelian group
$\integers^{d}$ with its standard generating set gives the graph of
Example~\ref{eg:lattice}, so we have that every metric boundary point
is a Busemann point.

\begin{figure}
    \begin{picture}(160,120)(-10,-60) \qbezier(0,0)(50,50)(100,50)
    \qbezier(0,0)(25,50)(100,50) \qbezier(0,0)(0,25)(100,50)
	
	\qbezier(0,0)(50,-50)(75,-50) \qbezier(0,0)(25,-50)(75,-50)
	\qbezier(0,0)(75,-25)(75,-50)

	\qbezier(100,50)(100,-50)(75,-50)
	\qbezier(100,50)(75,0)(75,-50)
	
	\put(0,0){\circle*{4}} \put(100,50){\circle*{4}}
	\put(75,-50){\circle*{4}}

	\put(-5,-5){\makebox(0,0)[tl]{$x$}}
	\put(105,55){\makebox(0,0)[br]{$y$}}
	\put(75,-55){\makebox(0,0)[t]{$z$}}
   \end{picture}
   \begin{picture}(160,120)(160,-60)
	
	\qbezier(200,0)(250,50)(300,50)
	\qbezier(200,0)(225,50)(300,50)
	
	\qbezier(200,0)(250,-50)(275,-50)
	\qbezier(200,0)(275,-25)(275,-50)

	\qbezier(300,50)(300,-30)(275,-50)
	
	\qbezier(225,0)(250,0)(300,50)
	\qbezier(225,0)(280,-10)(275,-50)
	
	\put(200,0){\line(1,0){25}}
	
	\put(200,0){\circle*{4}}
	\put(225,0){\circle*{4}}
	\put(300,50){\circle*{4}}
	\put(275,-50){\circle*{4}}

	\put(195,-5){\makebox(0,0)[tl]{$x$}}
	\put(305,55){\makebox(0,0)[br]{$y$}}
	\put(275,-55){\makebox(0,0)[t]{$z$}}
\end{picture}
    
    \caption{\label{fig:rigidtriple} Rigid and Non-rigid Triples}
\end{figure}
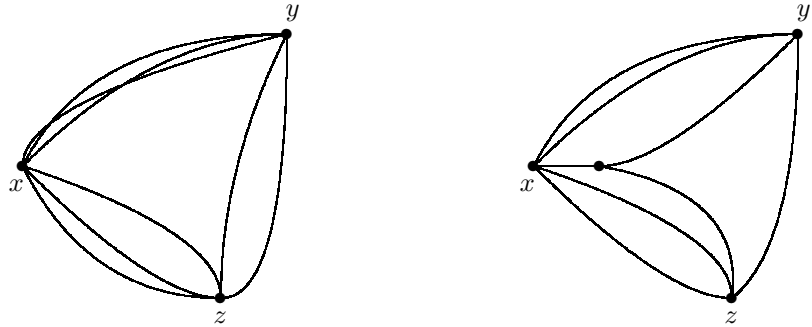

Since Cayley graphs have a lot more symmetry than arbitrary graphs, we
introduce the following definition which is easier to verify than the
condition of the previous section.

\begin{definition}
    Let $\Gamma = (V, E)$ be a graph.  A \emph{rigid triple} is a set
    of 3 vertices such that if $x$ is any one of these vertices, and
    $y$ and $z$ are the other two, then there are no minimal paths
    $[x,y]$ and $[x,z]$ which share a vertex other than $x$.
\end{definition}

Figure~\ref{fig:rigidtriple} illustrates the intuitive difference
between rigid and non-rigid triples.

In the Example~\ref{eg:nonbusemannpoints}, the sets $\{(1,1),
\gamma_{0}(n), (-1,1)\}$ are rigid triples, with one side having length
$2$, and the perimeter of these triples is $2n+2$, so we have an
infinite family of rigid triples of increasing perimeter but with one
side of bounded length.  The following two propositions show that this
sort of situation cannot occur in some important cases.

The utility of rigid triples is that if the only way they can have
large perimeter is if all the sides get large, we are guaranteed that
all boundary points are Busemann points.

\begin{proposition}\label{prop:pathsjoin}
    Let $\Gamma = (V,E)$ be a graph where for every $n$ there is an 
    $M_{n}$ such that every rigid triple $\{x,y,z\}$ with $d(x,y) = 
    n$ has perimeter bounded above by $M_{n}$, then if $\{a, b, c\}$ 
    is any triple of points with $d(a,b) = n$, and
    \[
      d(a,c) + d(b,c) > n + \max \{M_{k} : k = 1,\ldots,n\},
    \]
    then there are minimal paths $[a,c]$ and $[b,c]$ which share a
    tail.
\end{proposition}

\begin{proof}
    Assume that there are no such minimal paths.
    
    Let $[a,b]$, $[a,c]$, $[b,c]$ be minimal paths which
    maximize the total number of edges in common between $[a,b]$ and
    $[a,c]$, and between $[b,a]$ and $[b,c]$.  Clearly this
    total number of edges is at most $n-1$.  Let $x$ be the last
    vertex in common between $[a,b]$ and $[a,c]$, and $y$ the last
    vertex in common between $[b,a]$ and $[b,c]$.
    
    It is clear that $\{x,y,c\}$ is a rigid triple, otherwise we could
    either find more common edges, or find minimal paths $[a,c]$ and
    $[b,c]$ which share a tail.  Let $d(x,y) = k \le n$.  The perimeter of
    the triple is
    \begin{align*}
      d(x,c) + d(y,c) + d(x,y) &= (d(a,c) - d(a,x)) + (d(b,c) -
      d(b,x)) + \\
      & \qquad(d(a,b) - d(a,x) - d(b,y)) \\
      &= (d(a,c) + d(b,c)) + d(a,b) - 2(d(a,x) + d(b,y)) \\
      &> M_{k} + n + n - 2(n-1) > M_{k}
    \end{align*}
    Hence $\{x,y,c\}$ cannot be a rigid triple, and we have a
    contradiction.
\end{proof}

\begin{corollary}\label{cor:boundedrigidtriples}
    Let $\Gamma = (V,E)$ be a graph where for every $n$ there is an
    $M_{n}$ such that every rigid triple $\{x,y,z\}$ with $d(x,y) = n$
    has perimeter bounded above by $M_{n}$, then every point on the
    metric boundary is a Busemann point.
\end{corollary}

For Cayley graphs of groups, this is again a precise characterization
of the situation where all boundary points are Busemann points.

\begin{proposition}
    Let $G$ be a finitely generated group, with generating set $S$. 
    If for the corresponding Cayley graph, there is some $n$ such
    that for all $m$, there is a rigid triple $\{x,y,z\}$ with $d(x,y)
    = n$ and perimeter greater than $m$, then there is a point on the
    metric boundary which is not a Busemann point.
\end{proposition}

\begin{proof}
    Without loss of generality, since this is a Cayley graph, for each
    $m$ we can find rigid triples $\{e,y_{m},z_{m}\}$ with $d(e,y_{m})
    = n$ such that the perimeter greater than $m$.  Since the valence
    of each vertex of the graph is finite, the ball of radius $n$
    contains a finite number of points, and so there must be some
    particular $y$ with $d(e,y) = n$ such that $y = y_{m}$ for
    infinitely many $m$, and the minimal paths $[e,z_{m}]$,
    $[y,z_{m}]$ clearly do not share a tail.  Thus $e$ and $y$ are a
    pair of vertices which satisfy the conditions of
    Theorem~\ref{thm:nonbusemann}.  Hence the metric boundary of the
    graph contains a point which is not a Busemann point.
\end{proof}

This proposition will also hold in the setting of an arbitrary graph of finite valence
with an automorphism group which acts transitively on the vertices.

We can use this last Proposition to immediately prove the following:

\begin{proposition}
    Let $G$ be a group presented by a finite set of generators $S =
    S^{-1}$ with $|S|= k$ such that $G \isom F_{k}/N$, where $N$ is
    the normal group of all words in $S$ which map to $e$.  Let $N$ be
    finitely generated, and $M$ the maximum length of a generator of
    $N$.
    
    If $\Gamma$ is the Cayley graph of $G$ corresponding to these
    generators, no rigid triple has perimeter greater than $3M/2$, and
    so every metric boundary point is a Busemann point.
\end{proposition}

\begin{proof}
    Assume that there is a rigid triple whose perimeter $P$ is greater
    than $3M/2$.  We first note that from the triangle inequality the
    distance between any two of the vertices is at most $P/2$, and
    that the distance between at least one pair of vertices must be
    greater than or equal to $P/3$.
    
    Without loss of generality, we may assume that one vertex is the
    identity $e$ of $G$.  Let $x$ and $y$ be the other two vertices,
    and we can assume that $d(e,x) \ge P/3$.
    
    Let $w_{x}$, $w_{x^{-1}y}$ and $w_{y^{-1}}$ be minimal words
    representing $x$, $x^{-1}y$ and $y$, so that $w =
    w_{x}w_{x^{-1}y}w_{y^{-1}}$ represents a perimeter of the rigid
    triple.  Since the triple is rigid, there can be no cancellation
    within the product, so $w$ is a reduced word in $F_{k}$, and so we
    have that $w \in N$.  But $l(w) > M$, so $w$ cannot be a generator
    of $N$, and so we can write $w = g_{1}g_{2}\ldots g_{n}$, where
    $g_{k} \in R$, and so $l(g_{k}) \le M$.
    
    Pushing the $g_{k}$ down to $G$, we can see that each $g_{k}$ is a
    loop in the Cayley graph which starts and ends at $e$.  Moreover,
    $x$ must lie on at least one of these loops, say the loop
    corresponding to $g_{k}$.  But since this loop has length at most
    $M$, $d(e,x) \le M/2$.
    
    So we have $d(e,x) \ge P/3 > M/2 \ge d(e,x)$, which is a
    contradiction.  Hence no rigid triple has perimeter greater than
    $3M/2$.
\end{proof}

Not every finitely generated group has a finitely generated group of
relations.  Nevertheless, in some cases even when the defining set of
relations is infinite, there is a limit on the size of rigid triples. 
For example, $\integers^{n}$ with the standard generators can be
easily seen to have no rigid triples.

As might be expected, having an absolute bound on the size of rigid
triples is an exceptional situation, even for very nice presentations
of groups.

\begin{example}
  Let $G$ be the subgroup of $\complex$ generated by $e^{k \pi i/3}$
  for $k = 0, 1, \ldots, 5$.  $G$ is isomorphic to $\integers^{2}$,
  but the Cayley graph consists of the vertices and edges of a
  tessellation of the plane by equilateral triangles of side length 1.
  
  Triples of the form $\{0, k, ke^{\pi i/3}\}$ are rigid for all $k$,
  so we have rigid triples of arbitrarily large perimeter.  However, 
  one can see that if the distance between two points is $n$, then the 
  maximum size of a rigid triple is $3n$, and so every boundary point 
  is a Busemann point.
\end{example}

The above example is a special case of a more general phenomenon.  For
finitely generated abelian groups $G$, we do know that $G \isom
\integers^{d}/K$, where $\integers^{d}$ is the free abelian group
given by the generators, and $K$ is the ideal of all relations of $G$. 
A theorem of Dedekind says that $K$ is itself a finitely generated
free abelian group.  In the example, we have that $G \isom
\integers^{3}/K$, where $K = \{ (\lambda,-\lambda,\lambda) : \lambda
\in \integers \} \isom \integers$.

Using these facts, we can prove geometrically a result originally
proved by Develin (\cite{develin:cayleycompactification}, Theorem 7)
using algebraic and combinatoric techniques.

\begin{proposition}
    If $G$ be a finitely generated abelian group, then every point on
    the metric boundary of $G$ is a Busemann point.
\end{proposition}

\begin{proof}
    From Dedekind's theorem, we know that we can find $d$ and $K
    \subset \integers^{d}$ with $K \isom \integers^{l}$ for some $l$,
    such that $G \isom \integers^{d}/K$
    
    Let $\{e,x,y\}$ be a rigid triple in $G$, with $d(e,y) = n$.  We
    know that we can lift $x$, $y$ and $e$ to points $x'$, $y'$ and
    $z'$ so that we have $d(e,x') = d(e,x)$, $d(x',y') = d(x,y)$ and
    $d(y',z') = d(y,e) = n$.  Because we are in $\integers^{d}$, we
    must have that $x'$ lies on a minimal path $[e,y']$, else we can
    find minimal paths $[e,x']$ and $[y',x']$ which share a tail, and
    the projection of these into $G$ gives minimal paths $[e,x]$ and
    $[y,x]$ which share a tail, contradicting the rigidity of our
    original triple.  Similarly, we must have that $y'$ lies on a
    minimal path $[x',z']$, and so we have a minimal path from $e$ to
    $z'$ passing through $x'$, and $y'$.  So $d(e,z') = d(e,x) +
    d(x,y) + d(y,e)$ equals the perimiter of the rigid triple.

    Now $z' \in K$, so there are some unique integers $\lambda_{j} \ge
    0$ such that $z' = \sum_{j=1}^{l} \lambda_{j} k_{j}$, where $k_{j}
    = (k_{j,1}, k_{j,2}, \ldots, k_{j,d})$ are a complete set of
    (linearly independent) generators of $K$.  Letting $v = z' - y'$,
    we have that $l(v) = n$.  But we also have that $y'$ cannot
    contain even one full multiple of any one of the $k_{j}$, since
    the only way this can happen is if $x' = x'' + k_{1}'$ and $y' -
    x' = k_{2}' + y''$, where $k_{j} = k_{1}' + k_{2}'$.  But then if
    $l(k_{1}) \ne l(k_{2})$ we have that one of $x'' - k_{2}$ or
    $-k_{1} + y''$ is a shorter path which projects onto $[e,x]$ or
    $[x,y]$ respectively, which contradicts the construction of $x'$
    or $y'$; and if they are equal, then $x' - y' = - y'' + k_{1}'$
    gives a path which projects onto a minimal path $[y,x]$ which
    shares a tail with $[e,x]$, which will contradict the rigidity of
    our original triple.
    
    Therefore any minimal path $[e,y]$ must be a finite sum of
    elements of the form $a_{i} - v_{i}$, where $a_{i} = k_{j}$ for
    some $j$, $v = \sum v_{i}$, and $l(v) = \sum l(v_{i})$.  So there
    can be at most $n$ elements $a_{i}$ and each $a_{i}$ has length at
    most $(\max_{j} l(k_{j})) - 1$, so the length of $[e,y']$ is at most $n
    ((\max_{j} l(k_{j})) - 1)$.
    
    But then
    \[
      l(z') = l(v) + l(y') \le n + n ((\max_{j} l(k_{j})) - 1) \le n
      (\max_{j} l(k_{j})).
    \] 
    But $l(z')$ equals the original perimeter of the rigid triple, and
    we have determined a simple bound for it in terms of $n$.  Hence
    by Corollary~\ref{cor:boundedrigidtriples}, there are no
    non-Busemann points on the boundary.
\end{proof}

We conclude with some examples.

\begin{example}
    The free product group $G = \integers_{2} \ast \integers_{3}$ is
    generated by $a$, $a^{-1}$, $b$ and $b^{-1}$, and the relations
    $a^{2} = e$, $b^{3} = e$.  The largest rigid triples are the 
    triples $\{a, ab, ab^{2}\}$, which have perimeter 3.  Hence by 
    Proposition~\ref{prop:pathsjoin}, every metric boundary point is 
    a Busemann point.
    
    Similar analysis shows that for any finite free product
    $\ast F_{k}$ of finite groups $F_{k}$, $k =
    1,\ldots,n$, with generators the disjoint union $S =
    \bigcup_{k=1}^{n} S_{k}$ where $S_{k}$ generates $F_{k}$, the
    Cayley graph has rigid triples of size at most $\max |F_{k}|$, and
    so every metric boundary point is a Busemann point.
\end{example}

\begin{example}
    The braid group on $n$ strands, $B_{n}$, is given by generators
    $\{\sigma_{k}, \sigma_{k}^{-1} : k = 1, \ldots, n-1\}$ and
    relations $\sigma_{k}\sigma_{k+1}\sigma_{k} =
    \sigma_{k+1}\sigma_{k}\sigma_{k+1}$ for $k = 1, \ldots, n-2$, and
    $\sigma_{j}\sigma_{k} = \sigma_{k}\sigma_{j}$ for $|j - k| > 1$.
    
    Considering $B_{3}$, we note that the relations imply that
    $\sigma_{1}\sigma_{2}\sigma_{1}^{-1} =
    \sigma_{2}^{-1}\sigma_{1}\sigma_{2}$, and so
    $\sigma_{1}\sigma_{2}^{n}\sigma_{1}^{-1} =
    \sigma_{2}^{-1}\sigma_{1}^{n}\sigma_{2}$.
  
    Exactly as in Example~\ref{eg:nonbusemangroup}, we have that the
    triples $\{\sigma_{1},\sigma_{2}^{-1},
    \sigma_{1}\sigma_{2}^{n}\sigma_{1}^{-1}\}$ satisfy the conditions
    of Theorem~\ref{thm:nonbusemann}, and so this braid group has a
    non-Busemann point in its metric boundary.  By the same idea, we
    see that the families of triples $\{\sigma_{1},\sigma_{2}^{-1},
    \sigma_{1}\sigma_{2}^{-n}\sigma_{1}^{-1}\}$,
    $\{\sigma_{1}^{-1},\sigma_{2},
    \sigma_{1}^{-1}\sigma_{2}^{n}\sigma_{1}\}$, and
    $\{\sigma_{1}^{-1},\sigma_{2},
    \sigma_{1}^{-1}\sigma_{2}^{-n}\sigma_{1}\}$ satisfy
    Theorem~\ref{thm:nonbusemann}, so we have at least 4 distinct sets
    of non-Busemann points.
  
    The same construction will work for any pair of generators
    $\sigma_{k}$ and $\sigma_{k+1}$ of a braid group $B_{n}$.
\end{example}

We comment here that if we add a generator $b =
\sigma_{1}\sigma_{2}\sigma_{1}^{-1}$, then just as in
Example~\ref{eg:lipschitzequivalent} we get a new metric which is
Lipschitz equivalent to the original metric, but for which this
particular non-Busemann point becomes a Busemann point.  Hence even in
the special setting of Cayley graphs, Lipschitz equivalence does not
preserve whether or not a boundary point is a Busemann or non-Busemann
point.

It is unclear whether or not the existence of non-Busemann points for
Cayley graphs is invariant under Lipschitz equivalence; or even under
change of finite generating set.  It would be interesting to know the
answer to this question, but we conjecture based on
Example~\ref{eg:lipschitzequivalent} that for some groups at least,
the existence of non-Busemann points on the boundary will depend upon
the generating set.

\begin{figure}
    \begin{picture}(270,210)(-110,-110)
	\multiput(0,0)(30,0){5}{
	  \begin{picture}(100,100)(0,0)
	  \multiput(0,0)(-18,-18){4}{
	      \put(0,-30){\line(0,1){30}}
	      \put(0,0){\line(0,1){30}}
	      \put(0,0){\line(1,0){30}}
	      \put(0,0){\line(-1,0){30}}
	      \put(0,0){\circle*{2}}
	      \put(0,30){\line(0,1){30}}
	      \put(0,30){\circle*{2}}
	      \put(0,60){\line(0,1){30}}
	      \put(0,60){\circle*{2}}
	    }
	  \end{picture}
	}
	\multiput(0,0)(30,0){4}{
	  \begin{picture}(100,100)(0,0)
	  \multiput(0,0)(-18,-18){3}{
	      \put(0,30){\line(2,-3){12}}
	    }
	  \multiput(0,0)(-18,-18){2}{
	      \put(0,60){\line(-1,-6){6}}
	    }
	  \put(12,12){\circle*{4}}
	  \end{picture}
	}
	\put(5,-2){\makebox(0,0)[tl]{$e$}}
	\put(-13,-20){\makebox(0,0)[tl]{$x$}}
	\put(107,10){\makebox(0,0)[tl]{$c_{j}$}}
	
	\put(170,10){\vector(1,0){20}}
	\put(180,5){\makebox(0,0)[tc]{$k$}}
	
	\put(170,10){\vector(0,1){20}}
	\put(175,20){\makebox(0,0)[cl]{$m$}}
	
	\put(170,10){\vector(-1,-1){12}}
	\put(165,0){\makebox(0,0)[tl]{$n$}}
    \end{picture}
    
    \caption{\label{fig:heisenberg} Cayley Graph of the Heisenberg 
    Group}
\end{figure}
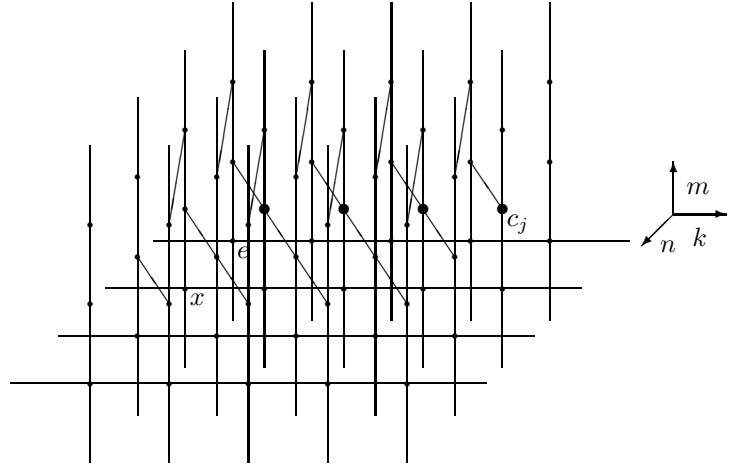

\begin{example}
  The discrete Heisenberg group is the following multiplicative
  subgroup of $GL_{3}$
  \[
    H^{3}_{d} = \left\{\begin{bmatrix}1 & m & n \\
    0 & 1 & k\\
    0 & 0 & 1\end{bmatrix} : m, n, k \in \integers\right\}.
  \]
  Let
  \[
    a = \begin{bmatrix}1 & 1 & 0 \\
    0 & 1 & 0\\
    0 & 0 & 1\end{bmatrix} \qquad \text{and} \qquad b =
    \begin{bmatrix}1 & 0 & 0 \\
    0 & 1 & 1\\
    0 & 0 & 1\end{bmatrix}.
  \]
  Then $S = \{a, b, a^{-1}, b^{-1}\}$ is a generating set for
  $H^{3}_{d}$.  The Cayley graph of this group with this generating
  set is illustrated in Figure~\ref{fig:heisenberg}.
  
  Triples $\{e, x = aba^{-1}b, c_{j} = b^{j-1}ab\}$ have
  perimeter $4 + 2(j+1)$, and there is a unique minimal path
  $b^{j-1}ab$ from $e$ to $b^{n-1}ab$, and a unique minimal path
  $b^{j}a$ from $aba^{-1}b$ to $b^{j-1}ab$.  These minimal paths share
  no common tail, so by Theorem~\ref{thm:nonbusemann}, there is a
  non-Busemann point on the metric boundary of the Cayley graph.
  
  Indeed, there are many non-Busemann points.  Triples of the form
  \[
    \{e, x, b^{j+1}a^{-1}b\}
  \]
  are also rigid and give distinct non-Busemann points from the above
  cases.  The group action by multiplication on the left give yet more
  examples, as do the triples formed by the inverses of these. 
  Another class of examples are triples of the form
  \[
    \{b, b^{-1}, b^{k}a^{-k}b^{-k}a^{k} =
  b^{-k}a^{k}b^{k}a^{-k}\},
  \]
  as well as left translates and inverses of these.
  
  There are some Busemann points, however.  The functions
  $\gamma^{v,\pm}_{n,k}: \naturals \to H_{d}^{3}$
  \[
    \gamma^{v,\pm}_{n,k}(t) = \begin{bmatrix} 1 & \pm t & n \\
      0 & 1 & k \\
      0 & 0 & 1 \\
    \end{bmatrix}
  \]
  are geodesic rays for fixed choice of $n$ and $k$, $+$ or $-$. 
  Similarly the functions $\gamma_{m,n,j}^{\pm}: \naturals \to
  H_{d}^{3}$ given by
  \[
    \gamma_{m,n}^{h,\pm}(t) = \begin{bmatrix} 1 & m & n \pm tm \\
      0 & 1 & \pm t \\
      0 & 0 & 1 \\
    \end{bmatrix}
  \]
  are geodesic rays for fixed choice of $n$, $k$ and $+$ or $-$.  All
  of these geodesic rays give distinct Busemann points.
\end{example}

The discrete Heisenberg group is of some significance, since it is
amenable, and hence has an amenable action on its metric boundary.  It
therefore may be susceptible to the sort of analysis that Rieffel used
in his discussion of $\integers^{d}$.  This would require finding
enough boundary points with finite orbits under the left action of the
group on the boundary.  Unfortunately, all the boundary points
described above have infinite orbits.

However, we have certainly not exhausted all possible boundary points
in this example.  For example, if we let $\omega^{v,+}_{n,k}$ be the
boundary point corresponding to the geodesic ray $\gamma^{v,+}_{n,k}$,
we conjecture that there is a boundary point or points of the form
\[
  \omega = \lim_{t \to \infty} \omega^{v,+}_{n_{t},k_{t}}
\]
where $n_{t} \to \infty$ and $k_{t} \to \infty$ as $t \to \infty$,
with $n_{t} \ge \alpha k_{t}$ eventually for any $\alpha$, and that
this point or points are fixed points of the action on the boundary.

Finally, it is common to consider
\[
  c = \begin{bmatrix}1 & 0 & 1 \\
    0 & 1 & 0\\
    0 & 0 & 1\end{bmatrix}
\]
as a generator of $H_{d}^{3}$ as well, and if this vertex is added to
$S$, many of the examples above remain rigid in the new metric.  It
may be that this metric is more amenable to study.  Of course 
we would like to be able to find the metric boundary for
arbitrary finite generating sets, but this problem seems difficult in
general.

\bibliographystyle{hplain}

\bibliography{busemann}

\end{document}